\documentclass[journal]{IEEEtran}

\usepackage[T1]{fontenc}
\usepackage[utf8]{inputenc}
\usepackage{color}
\usepackage[brazilian,english]{babel}

\usepackage{amsmath}
\usepackage{amssymb}
\usepackage{amsfonts}

\usepackage{graphicx}
\usepackage{psfrag}
\usepackage{subcaption}

\usepackage{amsthm}
\usepackage{thmtools}
\usepackage{thm-restate}

\usepackage{tikz}
\usepackage{pgfplots}
\usepgfplotslibrary{patchplots}
\usetikzlibrary{fit,positioning}

\usepackage{algpseudocode}

\usepackage{xspace}
\usepackage{multirow}
\usepackage{mathtools}
\usepackage{cancel}

\newcommand{\MCI}{\mathcal I}
\newcommand{\MCN}{\mathcal N}
\newcommand{\MCU}{\mathcal U}
\newcommand{\R}{\mathbb{R}}
\newcommand{\scale}{0.48}

\begin{document}

\title{Efficient Simultaneous Calibration of a Magnetometer and an
Accelerometer}
\author{Conrado S. Miranda,
Janito V. Ferreira
\thanks{C. S. Miranda is with the School of Electrical and Computer
  Engineering, University of Campinas, Brazil. E-mail:
conrado@dca.fee.unicamp.br.}
\thanks{J. V. Ferreira is with the School of Mechanical
Engineering, University of Campinas, Brazil. E-mail: janito@fem.unicamp.br.}}

\maketitle

\begin{abstract}
  This paper describes a calibration algorithm to simultaneously calibrate a
  magnetometer and an accelerometer without any information besides the sensors
  readings.
  Using a linear sensor model and maximum likelihood cost, the algorithm is able
  to estimate both sensors' biases, gains, and covariances, besides sensor
  orientations and Earth's fields.
  Results show errors of less than 0.1 standard deviations in simulation, and
  high-quality estimates with real sensors even when the algorithm's assumptions
  are violated.
\end{abstract}


\IEEEpeerreviewmaketitle

\section{Introduction}
\IEEEPARstart{S}{ensor} reading is an important step in a robot's control loop.
However, the values obtained may be inaccurate because of incorrect sensor
calibration. In the field of aerial vehicles, accelerometers and magnetometers
are frequently used to estimate attitude~\cite{Lizarragac}. As both these types
of sensors are calibrated using Earth's gravitational or magnetic field as
reference, most research into calibration of these devices focuses on
magnetometer calibration, which can also be used to calibrate accelerometers in
stationary settings. This work focuses on batch algorithms, that is, algorithms
that collect a data set and then calibrate the sensors, which can be used to
provide the initial conditions for online methods that solve the calibration
problem by filtering~\cite{Sabatini2006} or
iterating~\cite{Pylvanainen2008,Crassidis2005}.

A common approach to calibration is to assume that the magnetic field's norm is
known or unitary, and then adjust the gain and bias so that the readings match
this value. \cite{Alonso2002} extended the TWOSTEP algorithm, which tries to
minimize the difference between the field's norm and the sensed norm, to
calibrate both the bias and gain. \cite{vasconcelos2011geometric} proved that
noiseless magnetometer readings, including soft- and hard-iron effects, occur on
the surface of an ellipsoid manifold, and then developed an algorithm to
calibrate the bias and gain.

Although both methods can be used to calibrate accelerometers as well as
magnetometers, the calibration must be performed independently, so the
magnetometer and accelerometer are calibrated in different frames, and the
rotation between them is unknown. Furthermore, neither method is able to
estimate the direction of the magnetic field, which may be unknown and is
essential for attitude estimation using magnetometers
\cite{Sabatini2006,Gebre-Egziabher,Madgwick2011}.

\cite{Kok2012} partially solves this issue by calibrating the magnetometer while
also capturing a reference for the $z$ direction, provided by some independent
filter. This allows the direction of the Earth's magnetic field to be estimated,
and the calibration to be performed in some previously known reference frame,
while also estimating the magnetometer's bias and gain. However, a reference
orientation, whether reliable or not, is not always available, thus limiting the
use of algorithm.

\cite{Miranda2013} tried to solve both problems by simultaneously calibrating
the magnetometer and accelerometer, allowing the magnetic field, the attitude,
and all the sensors' parameters to be estimated. This eliminates the need for a
reference so that the calibration can be performed from sensor readings alone,
that is, no external apparatus or knowledge besides the sensor readings is
required. However, the problem is replete with low-quality local minima, which
sometimes significantly degrade the calibration performance, and the solution
has a very high computational cost, making the algorithm impractically slow.

The novelties in this paper are the improvement over \cite{Miranda2013} made
possible by the use of better initial estimates of the parameters, leading to
better minima, and the approximation for the rotation estimation, which also
helps to avoid poor minima and reduce computational time because of its
closed-form solution, making it suitable for use in practice. The initial
conditions are given by an adaptation of \cite{Kok2012}, where the $z$ reference
is replaced by an estimate computed from the accelerometer. To our knowledge,
the algorithms described in the present work and in \cite{Miranda2013} are the
only accelerometer and magnetometer calibration algorithms that are able to
perform the complete calibration from collected data alone. Moreover, the
algorithm proposed here is also the most statistically correct and generic
solution based on linear models since all the parameters (gains, biases, Earth's
fields, and orientations) must be estimated at the same time.

The algorithm is tested in simulated and real sensors. The simulations show that
the algorithm correctly calibrates different sets of sensor parameters and
allows comparison between the estimates and the ground truth. Several variations
of the algorithm are tested, and the best-performing one is also used to
calibrate the real sensors. The data capture was designed to violate the basic
assumptions underlying the algorithm so that its robustness to adverse
real-world conditions could be tested. It is shown that even in this imperfect
setting the proposed algorithm is able to achieve high-quality calibration.

The paper is organized as follows. Sections~\ref{sec:sensor_model}
and~\ref{sec:data_preprocessing} describe the sensor model used and the
collection and preprocessing of the data, respectively. Based on the sensor
model, a maximum likelihood cost function is described in
Section~\ref{sec:cost_function}. The initial parameters' estimates are given in
Section~\ref{sec:initial_estimate}, and the optimization steps are described in
Section~\ref{sec:optimization_steps}. The results of simulated and real-world
experiments are reported in Section~\ref{sec:experiments}  and are used to
evaluate the proposed method. The concluding remarks are presented in
Section~\ref{sec:conclusion}.

\section{Sensor Model}
\label{sec:sensor_model}
The linear sensor model with Gaussian noise, where the correct value is scaled
by a matrix and added to a bias and a noise, is frequently used in the
literature~\cite{Alonso2002,vasconcelos2011geometric} and is also used in this
paper. The value $v_s \in \R^3$ read by the sensor $s$, where $s=a$ and $s=m$
for the accelerometer and magnetometer, respectively, is given by:
\begin{equation}
  \label{eq:affine_measurements}
  v_s = K_s v_s^* + b_s + \epsilon_s, \quad \epsilon_s \sim \MCN(0, \Sigma_s),
\end{equation}
where $K_s \in \R^{3 \times 3}$ and $b_s \in \R^3$ are the gain matrix and bias
vector associated with the sensor, $v_s^* \in \R^3$ is the nominal value to be
read, and $\Sigma_s \in \mathcal S^3_{++}$ is the noise covariance, which is a
positive-definite matrix. It is important to note that the linear sensor model
for the magnetometer is also able to handle soft- and hard-iron effects, which
can be embedded in the sensor's gain and bias~\cite{vasconcelos2011geometric}.

Conversely, given the parameters and a value read by the sensor, the nominal
value can be estimated as:
\begin{equation*}
  v_s^* = K_s^{-1} (v_s - b_s) + \epsilon_s^*, \quad
  \epsilon_s^* \sim \MCN(0, K_s^{-1} \Sigma_s K_s^{-T}),
\end{equation*}
where it is assumed that the gain matrix $K_s$ is invertible. If it is not, then
its columns are not linear independent and the measured values $v_s$ must lie in
a reduced subspace spanned by these columns. Therefore, only the case where
$K_s$ is invertible is considered, and the nominal measurement $v_s^*$ can be
recovered.

The nominal value for a sensor is given by the Earth's respective field rotated
by the current orientation of the sensor frame, so that $v_s^* = R v_s^\MCI$,
where $v_s^\MCI$ is the field value for a given sensor in the fixed inertial
frame $\MCI$ and $R \in SO(3)$ is the rotation representing the current sensor
orientation. Since the choice of $\MCI$ is arbitrary, one can choose it such
that the gravity $g$ only has a negative component in $\vec z$ and the magnetic
field $h$ has a positive component in $\vec x$, such that
\begin{subequations}
\label{eq:fields}
\begin{gather}
  g = \begin{bmatrix}
    0 & 0 & g_z
  \end{bmatrix}^T,
  \quad
  h = \begin{bmatrix}
    h_x & 0 & h_z
  \end{bmatrix}^T,
  \\
  g_z < 0, \quad  h_x > 0,
  \\
  s_a^\MCI = g, \quad s_m^\MCI = h.
\end{gather}
\end{subequations}

It should be noted that the assumption that the magnetic field component $h_x$
is greater than zero can be satisfied almost everywhere except at the magnetic
poles, where the magnetic field is perpendicular to the Earth's surface and the
magnetometer does not provide more orientation information than the
accelerometer does for the gravitational field. In this particular case, one
cannot use this combination to perform the calibration, but one can do so in all
other contexts, as will be shown.

As discussed in~\cite{Miranda2013} and highlighted later, the sensor frame also
has multiple possible definitions, since the gain matrix can be written as $K_s
= R K_s'$, where $R$ is a rotation matrix. Therefore, some constraints must be
fixed to make the solution unique. This paper assumes that $K_a$ is triangular,
which is a sufficient assumption to define a single solution. Moreover, $K_m$ is
allowed to have any value, so that it can incorporate the misalignment between
the accelerometer and magnetometer.

Since during the calibration the algorithm does not use information about the
body frame, and this may not be available if the calibration is being performed
outside the desired body, other algorithms such as~\cite{Miranda2014} can be
used to align the sensor and body frame.

\section{Data Preprocessing}
\label{sec:data_preprocessing}
Assuming the linear model introduced in Sec.~\ref{sec:sensor_model} and that the
sensors can be held relatively still when readings are being taken so that a set
of sensor measurements have the same nominal value, the measured values can be
written as:
\begin{subequations}
\label{eq:sensor_reading_model}
\begin{gather}
  v_s[i,j] \sim \MCN(\mu_s[i], \Sigma_s),\quad
  \mu_s[i] = K_s R_i v_s^\MCI + b_s,
  \\
  j \in \{1,2,\ldots,\Delta_s[i]\},\quad
  i \in \{1,2,\ldots,N\},
\end{gather}
\end{subequations}
where $N$ is the number of measurement sets, $\Delta_s[i]$ is the number of data
points collected for the $i$-th set, $\mu_s[i]$ is the mean value shared by the
measurements in the $i$-th set, and $R_i$ is the orientation at which the $i$-th
set was collected.

Since the measurements in each set are assumed to have sampled the same
information, the mean values $\mu_s[i]$, which are different for each set, and
the covariance matrix $\Sigma_s$, which is shared among all sets, can be
estimated.

The minimum variance unbiased estimator for $\mu_s[i]$ is the sample
mean \cite{krishnamoorthy2006handbook}, given by:
\begin{equation}
  \label{eq:mu_hat}
  \hat \mu_s[i] = \frac{1}{\Delta_s[i]} \sum_{j=1}^{\Delta_s[i]} v_s[i,j],
\end{equation}
and has a distribution described by:
\begin{equation}
\label{eq:sensor_mean_distribution}
  \hat \mu_s[i] \sim \MCN \left(\mu_s[i], \frac{\Sigma_s}{\Delta_s[i]} \right).
\end{equation}
Since the sample mean $\hat \mu_s[i]$ characterizes the values for the $i$-th
set, it can be used instead of them as a single measurement with appropriate
weight, associated with the number of measurements $\Delta_s[i]$, during
optimization to decrease the computing time.

The minimum variance unbiased estimator for the noise covariance is given by:
\begin{equation}
\label{eq:covariance_initial_estimation}
  \hat \Sigma_s = \frac{\sum_{i=1}^N \sum_{j=1}^{\Delta_s[i]}
    \left(v_s[i,j] - \hat \mu_s[i] \right)
    \left(v_s[i,j] - \hat \mu_s[i] \right)^T
  }{\left(\sum_{i=1}^N \Delta_s[i] \right) - N},
\end{equation}
where the subtraction of $N$ from the denominator is a generalization of the
Bessel's correction \cite{kenney1957mathematics} for $N$ measurement sets,
instead of 1 as in the original correction.

\section{Cost Function}
\label{sec:cost_function}
Given that it is assumed that only the measured values are known, the full
parameter set is given by all the sensors' and fields' parameters and all
orientations. This set can be written as:
\begin{equation*}
  \Theta = \{\Sigma_a, K_a, b_a, g_z, \Sigma_m, K_m, b_m, h_x, h_z,
  R_1, R_2, \ldots, R_N \}.
\end{equation*}
Note that this is the most generic setting possible for this problem with the
assumption of affine measurements, since all parameters are considered unknown
and any previously known parameter can be fixed to the known value.

Using the model defined in Eq.~\eqref{eq:sensor_reading_model}, the
negative log-likelihood can be used to create a cost
function~\cite{bishop2006pattern}, which is given by:
\begin{subequations}
\label{eq:complete_cost}
\begin{gather}
  \begin{split}
    &J(\Theta; s_a, s_m) =
    \\
    &\quad
    \sum_{s \in \{a,m\}}
    \sum_{i=1}^N
    \sum_{j=1}^{\Delta_s[i]}
    \left(
      \ln(|\Sigma_s|) +
    \tilde v_s[i,j]^T \Sigma_s^{-1} \tilde v_s[i,j]
    \right)
  \end{split}
  \\
  \tilde v_s[i,j] = K_s R_i v_s^\MCI + b_s - v_s[i,j].
\end{gather}
\end{subequations}
In this format, the calibration problem becomes a traditional minimization
problem over all the parameters in $\Theta$, where the restrictions over
$\Sigma_s$, $g_z$, $h_x$, and $R_i$ must be satisfied.

A simplified cost function can be defined on a smaller parameter set $\Theta' =
\Theta \setminus \{\Sigma_a, \Sigma_m\}$, where the covariances are replaced by
their estimates computed using Eq.~\eqref{eq:covariance_initial_estimation}.
In this case, the maximum likelihood cost can be simplified to
\begin{subequations}
\label{eq:partial_cost}
\begin{gather}
  J'(\Theta'; \hat \mu_a, \hat \mu_m, \hat \Sigma_a, \hat \Sigma_m) =
  \sum_{i=1}^N \sum_{s \in \{a,m\}} \Delta_s[i]
  \tilde \mu_s[i]^T
  \hat \Sigma_s^{-1}
  \tilde \mu_s[i],
  \\
  \tilde \mu_s[i] = K_s R_i v_s^\MCI + b_s - \hat \mu_s[i],
\end{gather}
\end{subequations}
in which the real measurements $v_s[i,j]$ can be replaced by their sample means
$\hat \mu_s[i]$ with weight $\Delta_s[i]$, which greatly reduces the
computational cost.

Since both cost functions are subject to many local minima, which can lead to
poor estimates as shown in~\cite{Miranda2013}, a good initial estimate must be
provided. Then, instead of optimizing all the parameters at once, it will be
shown that optimizing subsets one at a time greatly simplifies the problem. The
next two sections deal with these problems.

\section{Initial Estimate}
\label{sec:initial_estimate}
This section is an adaptation of the work presented in~\cite{Kok2012}, where a
magnetometer was calibrated using a known correct reference measure of the $z$
direction. Since it is assumed that no such $z$ is available, this
reference is replaced here by the estimated accelerometer's field value, which
only has a negative component in the $z$ direction.

Furthermore, the constraint that the magnetometer cannot have mirrored axis
relative to the inertial axis is relaxed, which was not possible in the original
reference. The reader is referred to the original paper for a more detailed
derivation of the equations and implementation details.

Assuming that the field $v_s^\MCI$ is unitary and that the noise $\epsilon_s$
can be neglected, the following approximation can be performed:
\begin{subequations}
\label{eq:eta_components}
\begin{align}
  \label{eq:abc}
  0 &\approx \hat \mu_s[i]^T A \hat \mu_s[i] + b^T \hat \mu_s[i] + c \\
  A &= K_s^{-T} K_s^{-1} \\
  b &= -2 b_s^T K_s^{-T} K_s^{-1} \\
  c &= b_s^T K_s^{-T} K_s^{-1} b_s - 1,
\end{align}
\end{subequations}
which follows directly from taking the norm of both sides of
Eq.~\eqref{eq:affine_measurements} and replacing the single measurements by
their mean given by Eq.~\eqref{eq:mu_hat}, which has smaller variance and thus
is a better estimate.

By separating the known terms in Eq.~\eqref{eq:abc} from the unknown terms and
considering all $N$ sample sets, the approximation can be written as:
\begin{align}
  \label{eq:J_approx}
  J \eta &=
  \begin{bmatrix}
    J_1 \\
    \vdots \\
    J_N
  \end{bmatrix} \eta
  \approx 0
  \\
  J_i &= \begin{bmatrix}
    \hat \mu_s[i]^T \otimes \hat \mu_s[i]^T & \hat \mu_s[i]^T & 1
  \end{bmatrix}
  \nonumber
  \\
  \eta &= \begin{bmatrix}
    \text{vec} A \\
    b \\
    c
  \end{bmatrix}
  \nonumber
\end{align}
where $\otimes$ denotes the Kronecker product and $\text{vec}$ denotes the
vectorization operator. Assuming that $A$ is symmetric, the nontrivial solution
$\hat \eta$ that minimizes the approximation error is given by the right
eigenvector corresponding to the smallest singular value of $J$~\cite{Kok2012}.

Since any $\alpha \hat \eta$ is also a solution to the approximation in
Eq.~\eqref{eq:J_approx}, one can
find the correct $\alpha$ by manipulating the elements in
Eq.~\eqref{eq:eta_components}, which provides the correct scale. Hence the
correct value is given by:
\begin{equation*}
  \alpha = \left(\frac{1}{4} \hat b^T \hat A^{-1} \hat b - \hat c
  \right)^{-1}.
\end{equation*}

From the definitions in Eq.~\eqref{eq:eta_components} and the estimates $\hat
\eta$ and $\alpha$, the gain and bias must satisfy
\begin{align}
  \label{eq:decomposition_problem}
  \hat K_s^{-T} \hat K_s^{-1} &= \alpha \hat A
  \\
  \hat b_s &= -\frac{1}{2} \hat A^{-1} \hat b,
  \nonumber
\end{align}
which gives the bias estimate for each sensor directly.

The gain $\hat K_s$ can not be uniquely determined, as any $\hat K_s R$ with $R
R^T = I_3$ is also a solution to Eq.~\eqref{eq:decomposition_problem}. If the
upper triangular solution obtained by the Cholesky decomposition is denoted by
$\tilde K_s$, the estimates can be written as:
\begin{equation}
  \label{eq:sensor_coupling}
  \hat K_m = \tilde K_m R, \quad \hat K_a = \tilde K_a
\end{equation}
where it is assumed that only $\hat K_m$ is subject to rotation as the rotation
between the two gains is fixed and there would be multiple solutions if both had
full unknown rotations.

So far, the solution described in~\cite{Kok2012}, which provides estimates for
the sensor bias and gain, has been used to compute the initial estimates. This
solution could be applied here since the accelerometer and magnetometer could be
considered decoupled problems. However, in Eq.~\eqref{eq:sensor_coupling}, $\hat
K_a$ can be estimated directly while $\hat K_m$ depends on the unknown rotation
$R$. Since only the magnetometer was calibrated in~\cite{Kok2012}, this was not
an issue as the ground-truth $z$ direction was available. Moreover, the magnetic
field component $h_z$ could also be estimated directly. To solve this problem,
the ground truth $z$ was replaced by the estimate of the nominal value $v_a^*$
from the accelerometer, since this is already calibrated.

Therefore, estimates for the rotation $R$ between the accelerometer and the
magnetometer and the component $h_z$ of the magnetic field can be found
simultaneously by solving the following optimization problem:
\begin{subequations}
\label{eq:optimize_hz}
\begin{align}
  \min_{R,h_z} &&& \frac{1}{2} \sum_{i=1}^N
  \left(\Delta_a[i] + \Delta_m[i] \right)
  \left\|h_z + z_a[i]^T R^T z_m[i] \right\|^2_2
  \\
  \text{s.t.} &&&
  z_s[i] = \tilde K_s^{-1} \left(\hat \mu_s[i] - \hat b_s \right)
  \\
  &&& R \in SO(3)
\end{align}
\end{subequations}
where $\Delta_a$ and $\Delta_m$ are used to weight the number of samples in each
data set. As stated earlier, this equation is similar to the one
in~\cite{Kok2012}, but the ground truth $z$ is replaced by the estimates
$z_a[i]$.

Test showed that, unlike in \cite{Kok2012}, this problem is subject to a few
low-quality local minima because of the replacement of precise knowledge of the
$z$ direction by the accelerometer estimate and that the initial condition $R =
I_3$ and $h_z = 0$ sometimes leads to poor estimates. Hence, an optimization
with random restarts was performed and the initial conditions were given by a
random rotation as described in \cite{Kuffner2004} and $h_z \sim \MCU([-1,1])$.
Tests showed that at most $100$ iterations were needed to find good initial
estimates.

Since it was assumed that the fields were unitary, it follows that
\begin{equation*}
  \hat h_x = \sqrt{1-\hat h_z^2}.
\end{equation*}
As discussed in \cite{Miranda2013}, having the fields' and gains' scales as
variables leads to multiple solutions, since the gain can be multiplied by a
constant while the field is divided by the same value. Therefore, in this paper
the value $g_z = -1$ and $h_x=1$ are adopted, and the magnetometer estimates
must be adjusted as:
\begin{equation*}
  \hat K_m' = \hat K_m \hat h_x,
  \quad
  \hat h_z' = \hat h_z / \hat h_x,
\end{equation*}
while the accelerometer does not have to be adjusted because it already has a
unitary field.

Finally, the rotation $R_i$ associated with each sample set is given by the
eigenvector corresponding to the minimal eigenvalue of $B_i$ \cite{Karney2005},
where $B_i$ is defined as:
\begin{subequations}
\label{eq:rotation_approximation}
\begin{align}
  B_i &= w_a A_{a,i}^T A_{a,i} + w_m A_{m,i}^T A_{m,i}
  \\
  A_{s,i} &= f(\hat s_{s,i}^* + \hat s_s^\MCI,
               \hat s_{s,i}^* - \hat s_s^\MCI)
  \\
  \hat s_{s,i}^* &= \hat K_s^{-1} \left(\hat \mu_s[i] - b_s \right), ~
  w_s = |\hat K_s^T \hat \Sigma_s^{-1} \hat K_s|^\frac{1}{3}
  \\
  f(x,y) &= \begin{bmatrix}
    0 & -y_1 & -y_2 & -y_3 \\
    y_1 & 0 & -x_3 & x_2 \\
    y_2 & x_3 & 0 & -x_1 \\
    y_3 & -x_2 & x_1 & 0 \\
  \end{bmatrix}
\end{align}
\end{subequations}
and the weights $w_a$ and $w_m$ introduce uncertainty weighting.

All initial estimates therefore have closed-form efficient solutions except for
$R$ and $h_z$. However, Eq.~\eqref{eq:optimize_hz} has few decision variables
and can be efficiently optimized by random restarts and gradient descent.

\section{Optimization Steps}
\label{sec:optimization_steps}
As the cost functions in Eqs.~\eqref{eq:complete_cost} and
\eqref{eq:partial_cost} have many joint parameters, the parameter set is split
into four disjoint subsets that can be solved efficiently one at a time. The
optimization starts with the estimates provided in
Sec.~\ref{sec:initial_estimate} and then iterates each subset optimization while
keeping the other parameters fixed until convergence is achieved. Convergence is
defined by $J_{k-1} - J_k < \gamma$, where $J_k$ is the cost after the $k$-th
iteration and $\gamma$ is the stop parameter. Since the true values for any
parameter in $\Theta$ are not available, the hat notation is dropped where
convenient to simplify the expressions.

This section is an extension of the optimization presented
in~\cite{Miranda2013}. The main differences are, firstly, the approximation for
the rotation estimation described in Sec.~\ref{sec:rotation_estimation}. This
uses Eq.~\eqref{eq:rotation_approximation} to compute an efficient approximation
to the real orientation during the optimization, which is not present
in~\cite{Miranda2013}, leading to a significant slow-down of the algorithm and
worse results. Secondly, the format of the gain matrix $K_a$ computed in
Sec.~\ref{sec:optimization_steps:gain} is assumed to be triangular here because
its initial condition is given by a Cholesky decomposition while
in~\cite{Miranda2013} it is assumed to be symmetric. The main difference between
the method proposed here and the one described in~\cite{Kok2012}, apart from the
fact that two sensors are used and no ground-truth is available in the former,
is that the optimization is performed in steps with efficient solutions, while
in~\cite{Kok2012} the initial estimate is computed and then all parameters are
optimized at the same time, preventing the inherent structure of the problem
from being exploited.
\subsection{Rotation Estimation}
\label{sec:rotation_estimation}
From Eq.~\eqref{eq:partial_cost}, it is clear that the cost for the $i$-th
rotation is given by:
\begin{gather}
  \label{eq:rotation_estimation}
  J'_{R_i} =
  \sum_{s \in \{a,m\}} \Delta_s[i]
  \tilde \mu_s[i]^T
  \Sigma_s^{-1}
  \tilde \mu_s[i],
  \\
  \tilde \mu_s[i] = K_s R_i s_s^\MCI + b_s - \hat \mu_s[i].
  \nonumber
\end{gather}

As a closed-form solution could not be found, the problem is optimized directly
through gradient descent using a quaternion to represent the rotation. However,
the approximation used to compute the initial estimates, which is given by the
eigenvector associated with the minimum eigenvalue of $B_i$ from
Eq.~\eqref{eq:rotation_approximation}, can be used to speed-up the initial steps
of the optimization, when the estimates are still far from the true values. The
performance of this approximation approximation is evaluated in the experimental
section.
\subsection{Bias and Field Estimation}
Since the sensor readings are assumed linear in the biases and fields and have
Gaussian noise, the problem can be posed as a generalized least squares (GLS)
problem \cite{kariya2004generalized} given by finding $\beta$ to approximate
\begin{equation*}
  y = X \beta + \epsilon, ~\text{E}[\epsilon|X] = 0,
  ~\text{Var}[\epsilon|X] = \Omega,
\end{equation*}
whose solution is given by:
\begin{equation*}
  \hat \beta = \left(X^T \Omega^{-1} X \right)^{-1} X^T \Omega^{-1} y.
\end{equation*}

Writing $\hat \mu_a[i]$ in the GLS format using
Eqs.~\eqref{eq:sensor_reading_model} and \eqref{eq:sensor_mean_distribution}
gives
\begin{gather*}
  \begin{bmatrix}
    \hat \mu_a[1] \\
    \vdots \\
    \hat \mu_a[N]
  \end{bmatrix}
  =
  \begin{bmatrix}
    K_a R_1 \vec z & I_3 \\
    \vdots & \vdots \\
    K_a R_N \vec z & I_3
  \end{bmatrix}
  \begin{bmatrix}
    g_z \\
    b_a
  \end{bmatrix}
  + \epsilon,
  \\
  \Omega = \text{diag} \left( \frac{\Sigma_a}{\Delta_a[1]}, \ldots,
  \frac{\Sigma_a}{\Delta_a[N]} \right),
\end{gather*}
where $\vec z$ is the unitary vector in the $z$ direction. As
Sec.~\ref{sec:initial_estimate} established that $g_z = -1$, the gain matrix is
adjusted to ensure that this equality is still valid. The estimate is therefore
adjusted so that $\hat K_a' = -\hat K_a \hat g_z$, where $\hat g_z$ is the
estimated gravity.

The magnetometer estimate $\hat \mu_m[i]$ can be written in a similar fashion,
allowing $\hat h_x$, $\hat h_z$ and $\hat b_m$ to be estimated. As $h_x=1$, it
follows that $\hat K_m' = \hat K_m \hat h_x$ and $\hat h_z' = \hat h_z / \hat
h_x$, which maintains the scale.
\subsection{Gain Estimation}
\label{sec:optimization_steps:gain}
The gains can also be written in the GLS format using
Eqs.~\eqref{eq:sensor_reading_model} and \eqref{eq:sensor_mean_distribution}.
Assuming that $K_a$ is upper triangular, which is valid for the initial estimate
since it comes from a Cholesky decomposition, one has that:
\begin{gather*}
  \begin{bmatrix}
    \hat \mu_a[1] - b_a \\
    \vdots \\
    \hat \mu_a[N] - b_a
  \end{bmatrix}
  =
  \begin{bmatrix}
    G_1 \\
    \vdots \\
    G_N
  \end{bmatrix}
  \begin{bmatrix}
    K_{a,11} \\
    K_{a,12} \\
    K_{a,13} \\
    K_{a,22} \\
    K_{a,23} \\
    K_{a,33}
  \end{bmatrix}
  + \epsilon,
  \\
  G_i =
  \begin{bmatrix}
    g_{R_i,1} & g_{R_i,2} & g_{R_i,3} & 0 & 0 & 0
    \\
    0 & 0 & 0 & g_{R_i,2} & g_{R_i,3} & 0
    \\
    0 & 0 & 0 & 0 & 0 & g_{R_i,3}
  \end{bmatrix},
  \\
  g_{R,j} = (R g)_j, ~
  \Omega = \text{diag} \left( \frac{\Sigma_a}{\Delta_a[1]}, \ldots,
  \frac{\Sigma_a}{\Delta_a[N]} \right),
\end{gather*}
where $(v)_i$ represents the $i$-th component of the vector $v$.

Again, a similar formulation can be used for the magnetometer, but the gain
$K_m$ is considered full, as it is the product of a rotation with an upper
triangular matrix.
\subsection{Covariance Estimation}
The maximum likelihood covariance estimator is similar that given in
Eq.~\eqref{eq:covariance_initial_estimation}, but the reference $\hat \mu_s[i]$
is replaced by an estimate of $\mu_s[i]$ using parameters in $\Theta$, and
Bessel's correction is not used, since it lowers the likelihood of the data in
order to provide an unbiased estimate. The new estimate is then given by:
\begin{gather}
  \hat \Sigma_s = \frac{\sum_{i=1}^N \sum_{j=1}^{\Delta_s[i]}
    \tilde s_s[i,j] \tilde s_s[i,j]^T
  }{\sum_{i=1}^N \Delta_s[i]},
  \label{eq:covariance_estimation}
  \\
  \tilde s_s[i,j] = s_s[i,j] - \left(K_s R_i s_s^\MCI + b_s \right).
  \nonumber
\end{gather}

\section{Experiments}
\label{sec:experiments}
Two experiments are performed to validate the proposed algorithm: calibration of
simulated sensors and calibration of real sensors. The simulation is used to
evaluate the effects of changing the stop condition, the number of sampled
intervals, and the algorithm hypothesis. The performance is analyzed in terms of
errors in relation to the ground-truth values, which are available in the
simulation, and the total time required for the simulation. Using a Monte Carlo
approach, many combinations of sensors can be simulated to evaluate the
robustness of the proposed algorithm.

The performance with test data, which is not used for training, can be evaluated
to test the generalization of the algorithm~\cite{bishop2006pattern}.

In the real experiment, a magnetometer and an accelerometer had their readings
measured while being held by hand inside a building, which slightly violates the
hypothesis that all samples measure the same value. It will be shown that the
algorithm is able to fit the sampled data correctly.

It is important to note that only the algorithm derived in this paper is used
for calibration since the authors do not know of any other work, apart from
their own previous contribution~\cite{Miranda2013} upon which this paper is
based, that describes an algorithm that can perform this kind of calibration
without external references. If an algorithm such as the one proposed
in~\cite{Kok2012} were used, additional information would have to be provided,
making comparison of the results impossible since the calibration would be
performed on different data. Thus, no comparison with any previous work is
performed.
\subsection{Simulated calibration}
\begin{table}[!b]
  \centering
  \caption{Experimental parameters}
  \begin{tabular}{|c|c|c|}
    \hline
     & Description & Value
    \\
    \hline
    $N$ & \# of measurement sets & $15$
    \\
    $\Delta[i]$ & \# of samples for each set $i$ &
    $\MCU(\{400,401,\ldots,600\})$
    \\
    $R_i$ & Rotation for set $i$ & See \cite{Kuffner2004}
    \\
    $\gamma$ & Stop condition & $10^{-4}$
    \\
    $g_z$ & Gravity's $z$ component & $\MCU([-1.5,-0.5])$
    \\
    $h_x$ & Magnetic field's $x$ component & $\MCU([0.5,1.5])$
    \\
    $h_z$ & Magnetic field's $z$ component & $\MCU([-1.5,1.5])$
    \\
    $b_{s,i}$ & Bias' component & $\MCU([-1,1])$
    \\
    $K_s$ & Gain matrix & See Eq.~\eqref{eq:experiment_gain}
    \\
    $\Sigma_s$ & Covariance matrix & See Eq.~\eqref{eq:experiment_covariance}
    \\
    & \# of Monte Carlo runs & $100$
    \\
    \hline
  \end{tabular}
  \label{tab:experiment}
\end{table}

\begin{figure*}[!t]
  \centering
  \begin{subfigure}{\scale\textwidth}
    \centering
    \psfrag{NCAR---}[l][l]{\scriptsize NCAR}
    \psfrag{FCAR---}[l][l]{\scriptsize FCAR}
    \psfrag{DCAR---}[l][l]{\scriptsize DCAR}
    \psfrag{NCDR---}[l][l]{\scriptsize NCDR}
    \psfrag{stop}[t][c]{\footnotesize $\log_{10} \gamma$}
    \psfrag{yaxis}[b][c]{\footnotesize $\log_{10} \delta$}
    \includegraphics[width=\textwidth]{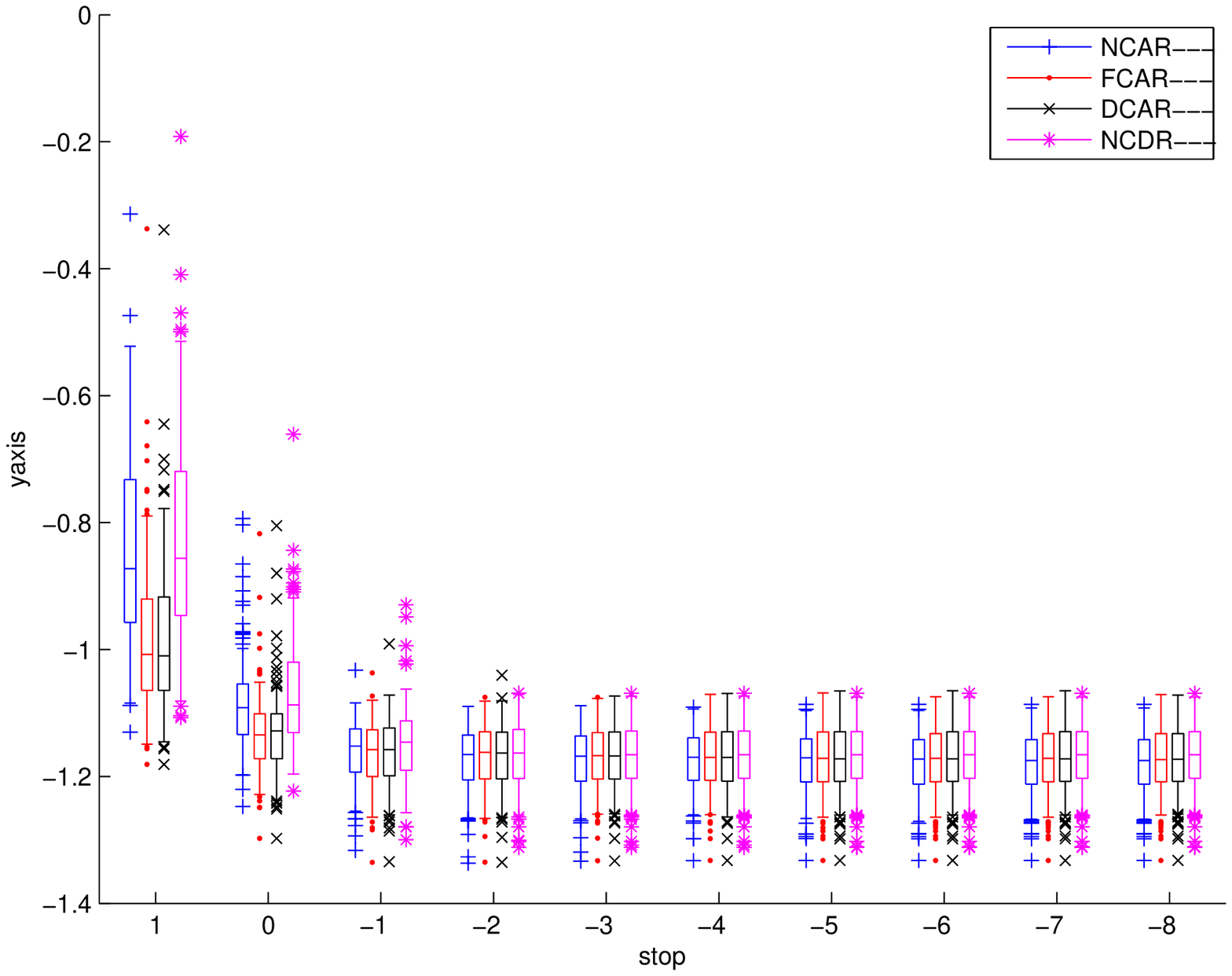}
    \caption{Accelerometer}
  \end{subfigure}
  ~
  \begin{subfigure}{\scale\textwidth}
    \centering
    \psfrag{NCAR---}[l][l]{\scriptsize NCAR}
    \psfrag{FCAR---}[l][l]{\scriptsize FCAR}
    \psfrag{DCAR---}[l][l]{\scriptsize DCAR}
    \psfrag{NCDR---}[l][l]{\scriptsize NCDR}
    \psfrag{stop}[t][c]{\footnotesize $\log_{10} \gamma$}
    \psfrag{yaxis}[b][c]{\footnotesize $\log_{10} \delta$}
    \includegraphics[width=\textwidth]{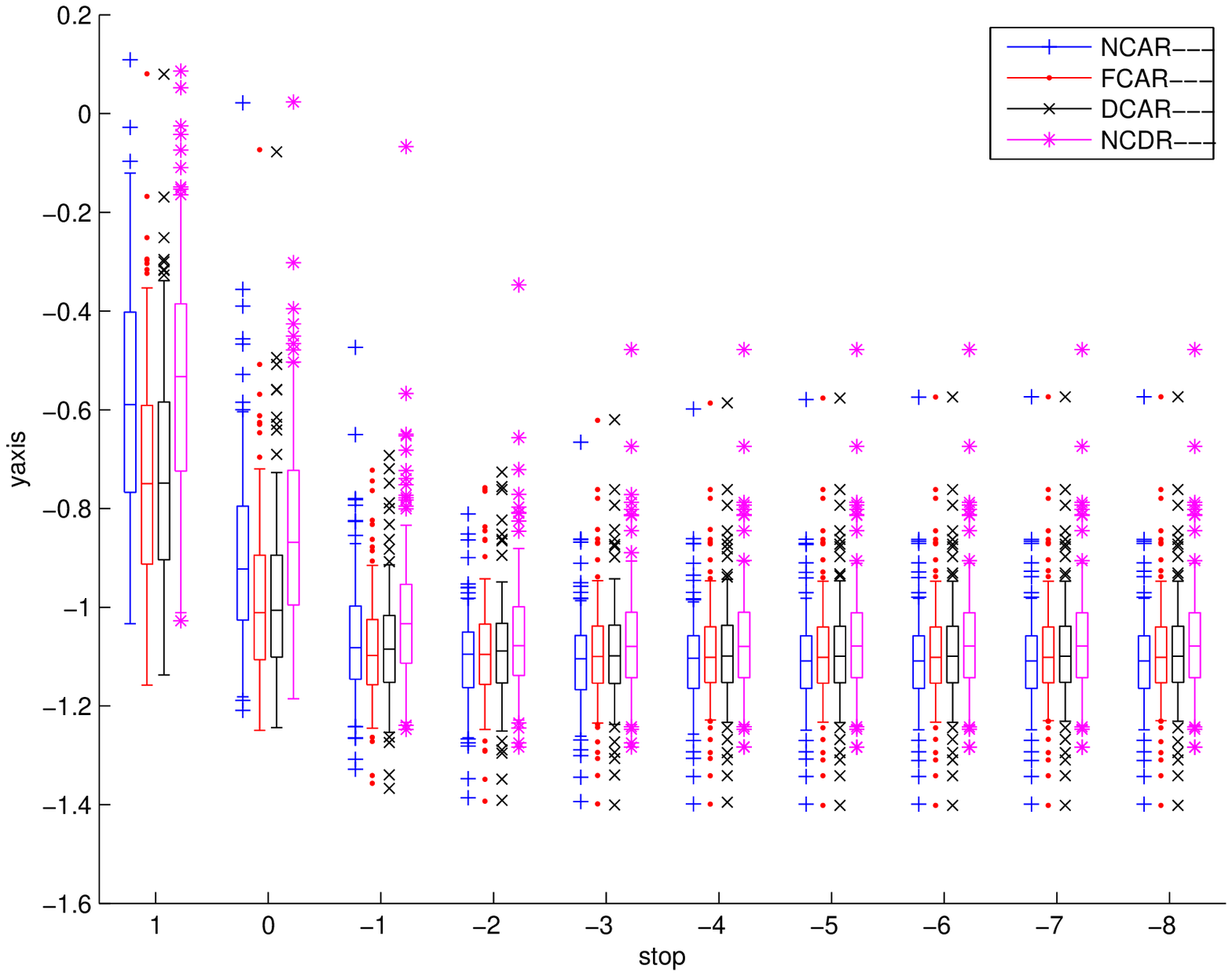}
    \caption{Magnetometer}
  \end{subfigure}
  \caption{Reconstruction error for training data for different stop conditions.
  Best viewed in color.}
  \label{fig:results:stop_train}
\end{figure*}

\begin{figure*}[!t]
  \centering
  \begin{subfigure}{\scale\textwidth}
    \centering
    \psfrag{NCAR---}[l][l]{\scriptsize NCAR}
    \psfrag{FCAR---}[l][l]{\scriptsize FCAR}
    \psfrag{DCAR---}[l][l]{\scriptsize DCAR}
    \psfrag{NCDR---}[l][l]{\scriptsize NCDR}
    \psfrag{stop}[t][c]{\footnotesize $\log_{10} \gamma$}
    \psfrag{yaxis}[b][c]{\footnotesize $\log_{10} \delta$}
    \includegraphics[width=\textwidth]{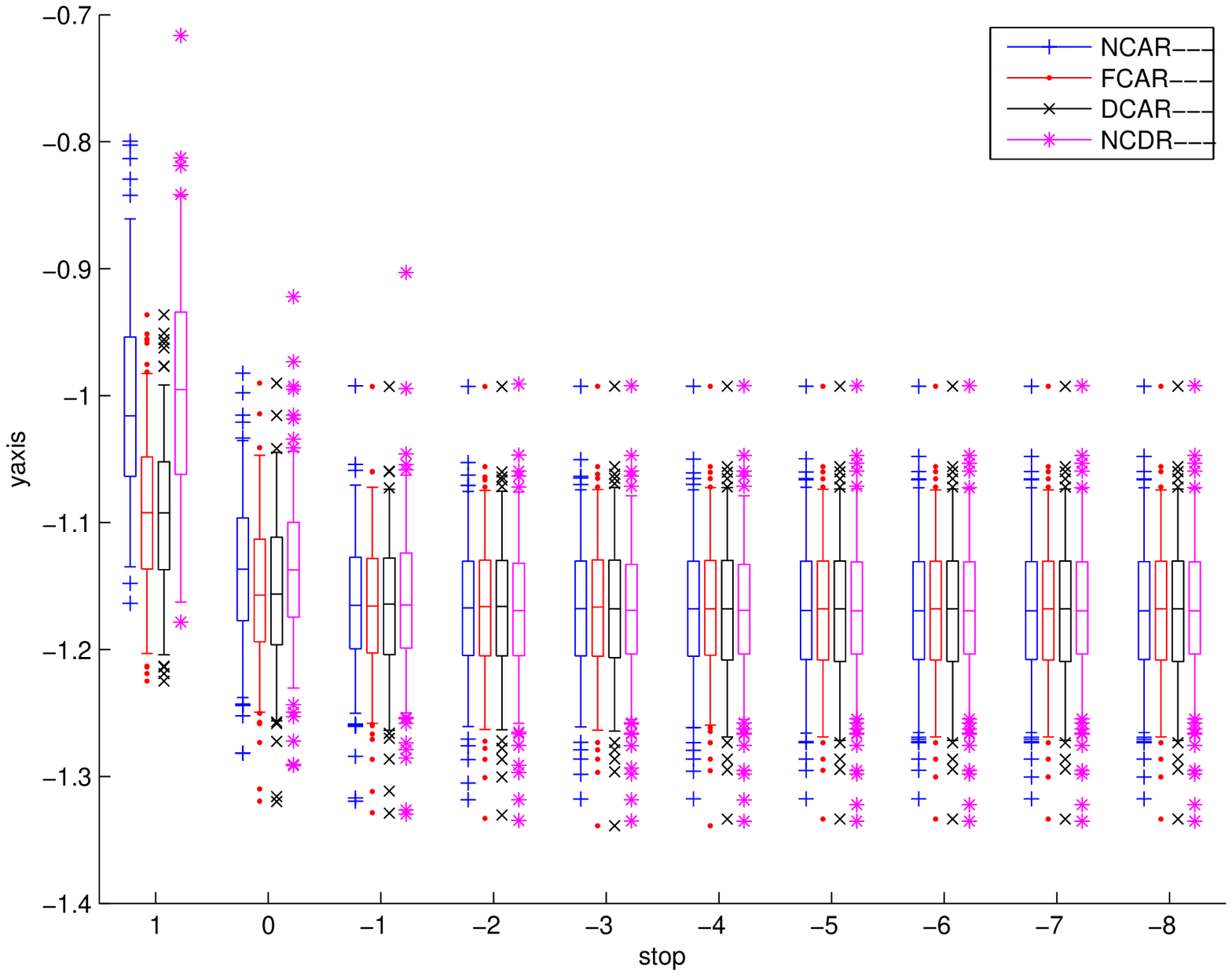}
    \caption{Accelerometer}
  \end{subfigure}
  ~
  \begin{subfigure}{\scale\textwidth}
    \centering
    \psfrag{NCAR---}[l][l]{\scriptsize NCAR}
    \psfrag{FCAR---}[l][l]{\scriptsize FCAR}
    \psfrag{DCAR---}[l][l]{\scriptsize DCAR}
    \psfrag{NCDR---}[l][l]{\scriptsize NCDR}
    \psfrag{stop}[t][c]{\footnotesize $\log_{10} \gamma$}
    \psfrag{yaxis}[b][c]{\footnotesize $\log_{10} \delta$}
    \includegraphics[width=\textwidth]{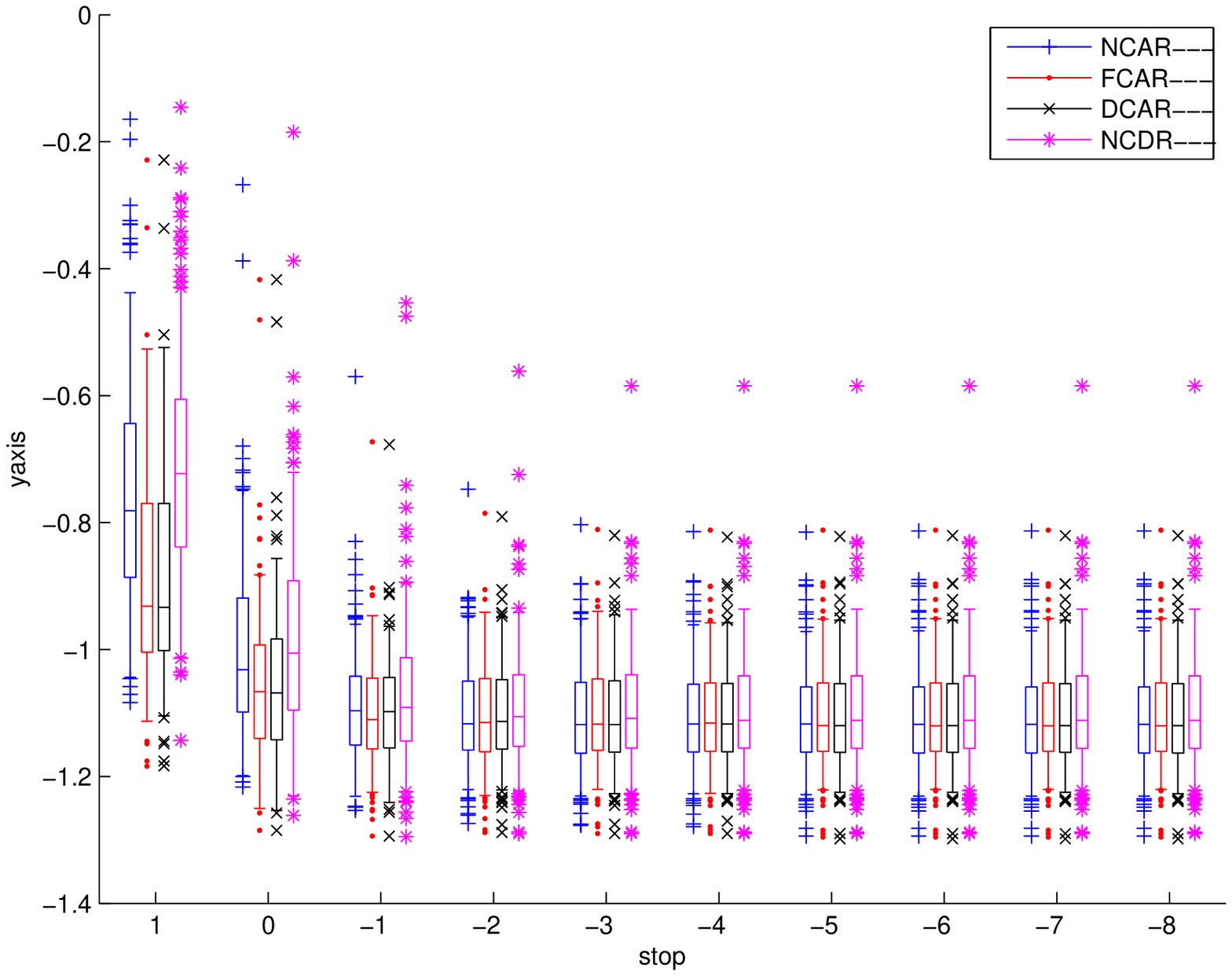}
    \caption{Magnetometer}
  \end{subfigure}
  \caption{Reconstruction error for test data for different stop conditions.
  Best viewed in color.}
  \label{fig:results:stop_test}
\end{figure*}

To evaluate the calibration algorithm, 100 Monte Carlo simulations with new
sensor parameters on each run were performed, allowing evaluation of both common
and possible exceptional calibration behavior. This analysis is important to
show that the method described handles different sensor characteristics well.
Table~\ref{tab:experiment} presents a summary of the parameters used in the
experiment; the symbols in the table are those used throughout the paper.

For each Monte Carlo run, a total of $N=15$ sets of sensor measurements were
considered, with the same number of samples $\Delta[i] \sim
\MCU(\{400,401,\ldots,600\})$ for both sensors. The rotations for each set were
randomly created using the method described in \cite{Kuffner2004}, and the stop
condition was set to $\gamma = 10^{-4}$.

When varying the stop condition $\gamma$, the same sample sets were used for
each value tested, so that the performance can be viewed as the same algorithm
stopping at different times. This avoids differences caused by random
fluctuations and allows better comparison. For a varying number of sets, runs
with $N+1$ sample sets merged a new set to the previous $N$ already used,
simulating collection of an increasing number of sets and also preventing random
fluctuations from significantly affecting the comparison. Alternatively, the
training with $N$ sets can be viewed as training on a subset of the data set
provided to the algorithm that trained with $N+1$ sets, thus maintaining
consistency.

\begin{figure*}[!t]
  \centering
  \begin{subfigure}{\scale\textwidth}
    \centering
    \psfrag{NCAR---}[l][l]{\scriptsize NCAR}
    \psfrag{FCAR---}[l][l]{\scriptsize FCAR}
    \psfrag{DCAR---}[l][l]{\scriptsize DCAR}
    \psfrag{NCDR---}[l][l]{\scriptsize NCDR}
    \psfrag{sample}[t][c]{\footnotesize $N$}
    \psfrag{yaxis}[b][c]{\footnotesize $\log_{10} \delta$}
    \includegraphics[width=\textwidth]{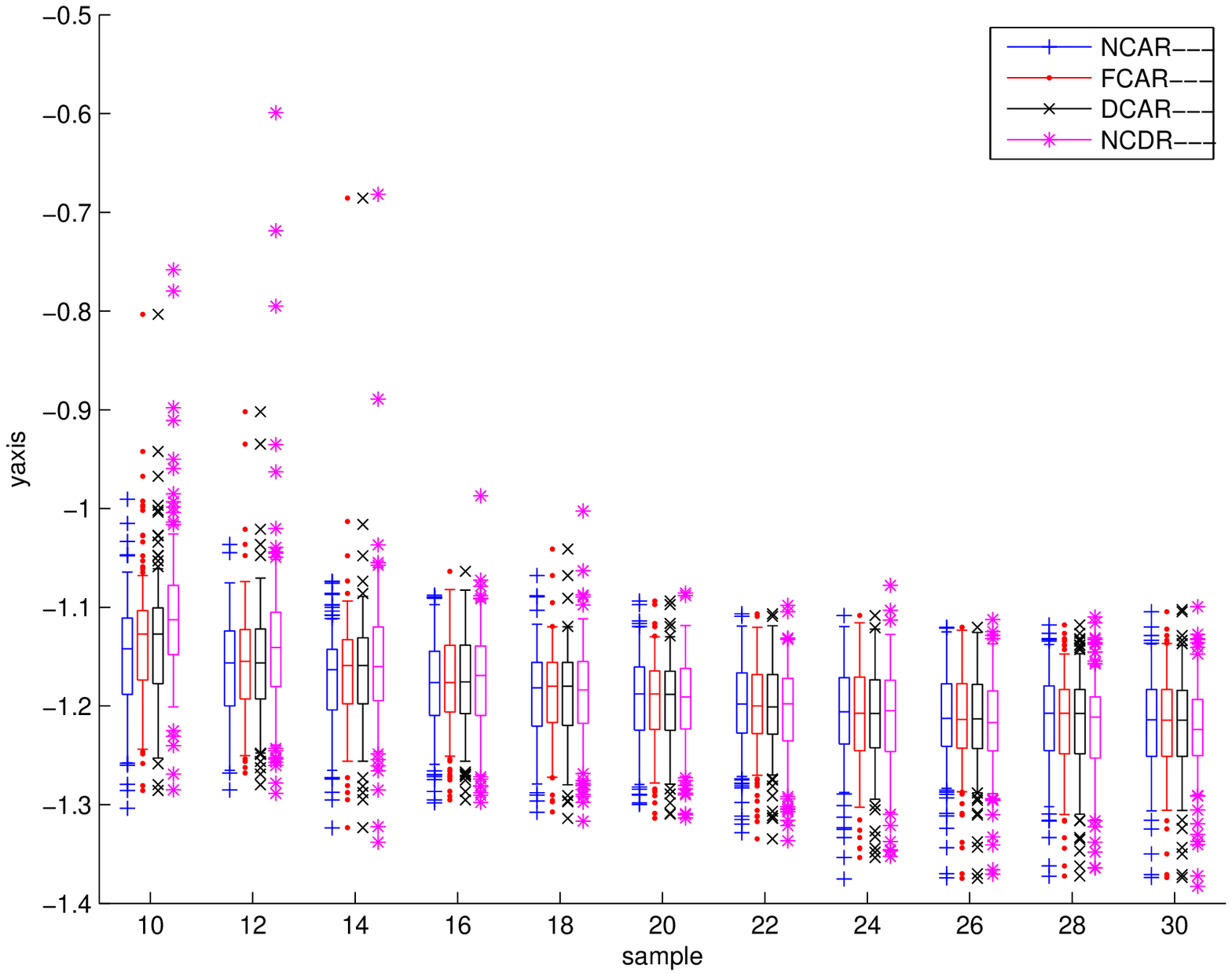}
    \caption{Accelerometer}
  \end{subfigure}
  ~
  \begin{subfigure}{\scale\textwidth}
    \centering
    \psfrag{NCAR---}[l][l]{\scriptsize NCAR}
    \psfrag{FCAR---}[l][l]{\scriptsize FCAR}
    \psfrag{DCAR---}[l][l]{\scriptsize DCAR}
    \psfrag{NCDR---}[l][l]{\scriptsize NCDR}
    \psfrag{sample}[t][c]{\footnotesize $N$}
    \psfrag{yaxis}[b][c]{\footnotesize $\log_{10} \delta$}
    \includegraphics[width=\textwidth]{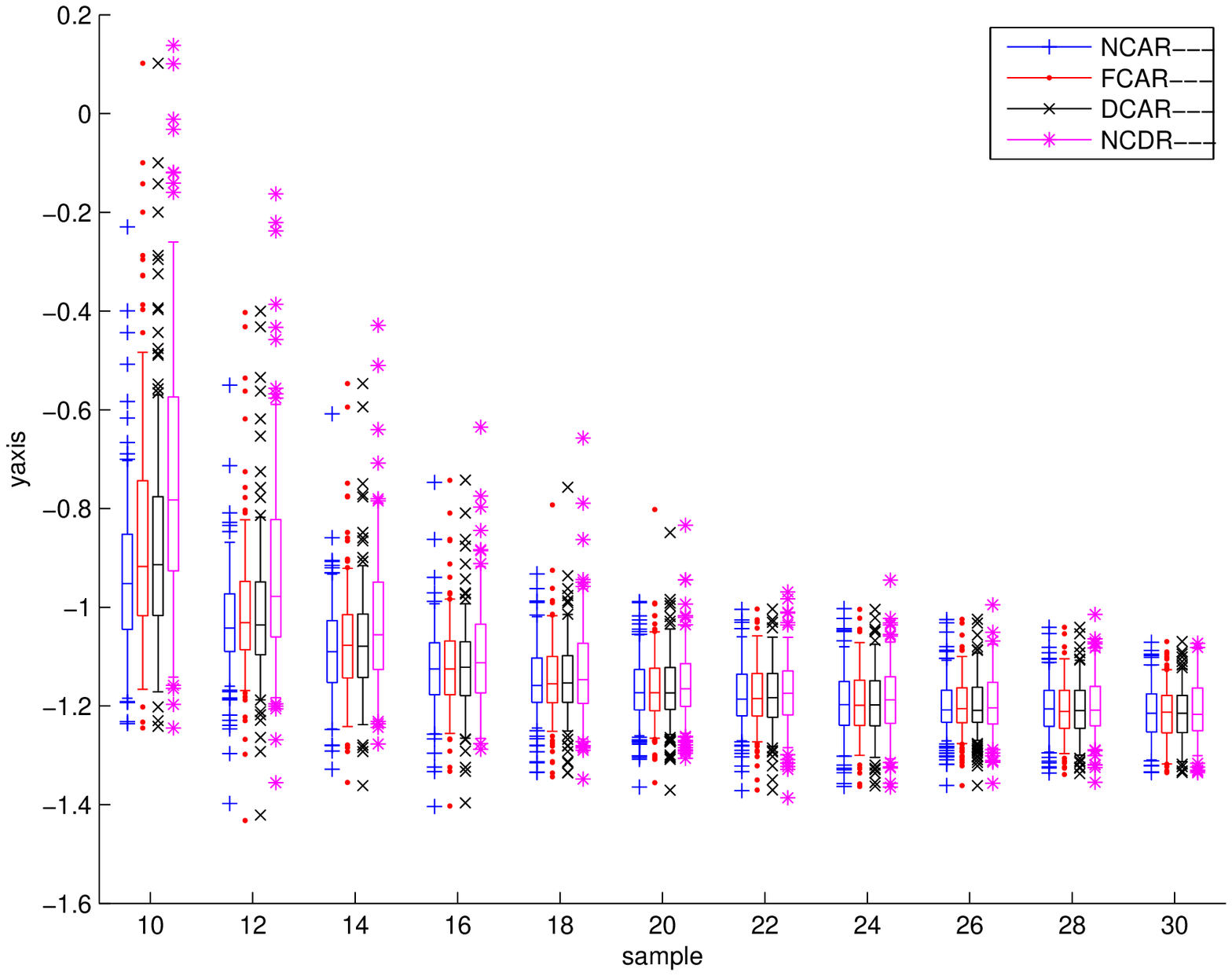}
    \caption{Magnetometer}
  \end{subfigure}
  \caption{Reconstruction error for training data for different numbers of
  sample sets. Best viewed in color.}
  \label{fig:results:sample_train}
\end{figure*}

\begin{figure*}[!t]
  \centering
  \begin{subfigure}{\scale\textwidth}
    \centering
    \psfrag{NCAR---}[l][l]{\scriptsize NCAR}
    \psfrag{FCAR---}[l][l]{\scriptsize FCAR}
    \psfrag{DCAR---}[l][l]{\scriptsize DCAR}
    \psfrag{NCDR---}[l][l]{\scriptsize NCDR}
    \psfrag{sample}[t][c]{\footnotesize $N$}
    \psfrag{yaxis}[b][c]{\footnotesize $\log_{10} \delta$}
    \includegraphics[width=\textwidth]{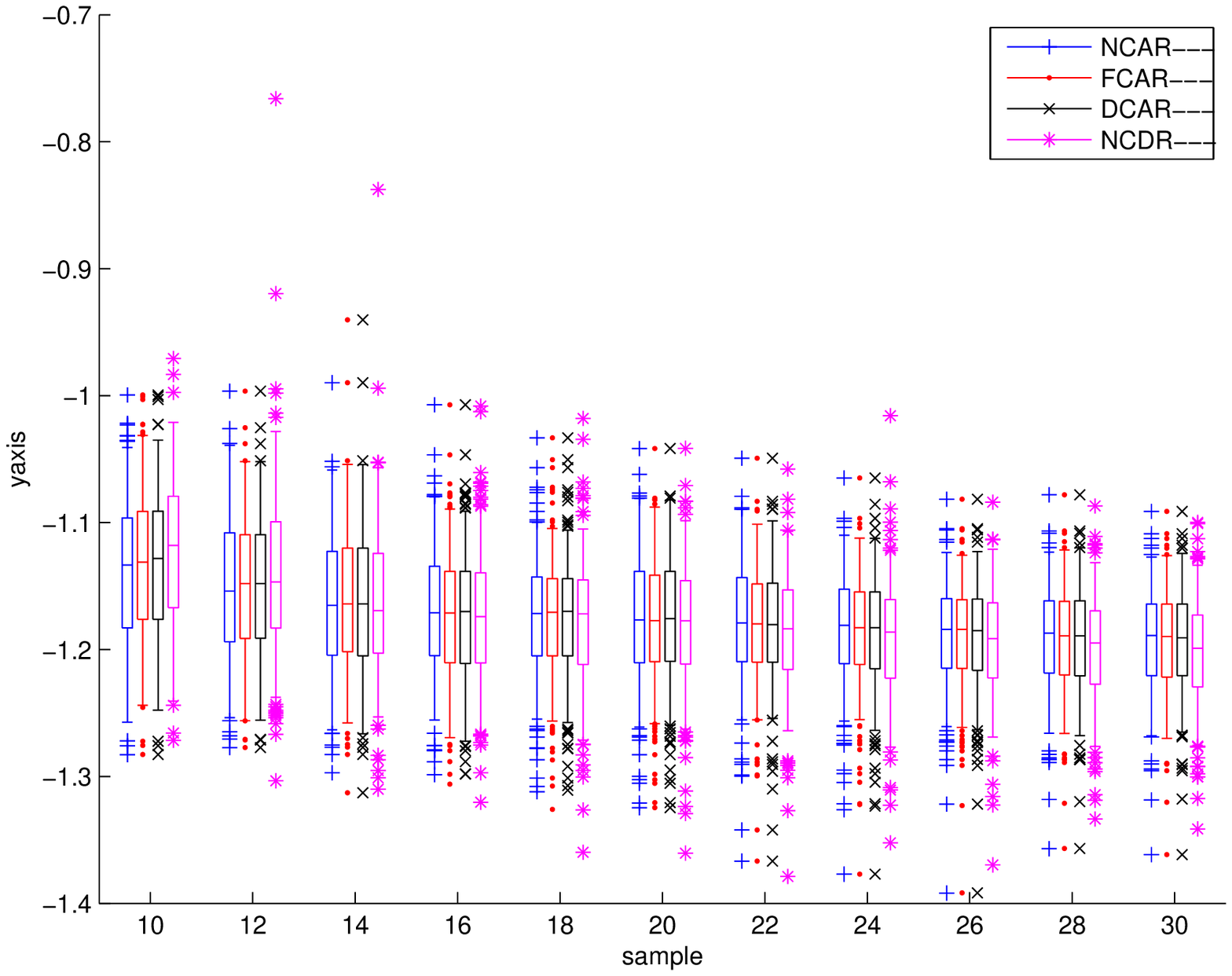}
    \caption{Accelerometer}
  \end{subfigure}
  ~
  \begin{subfigure}{\scale\textwidth}
    \centering
    \psfrag{NCAR---}[l][l]{\scriptsize NCAR}
    \psfrag{FCAR---}[l][l]{\scriptsize FCAR}
    \psfrag{DCAR---}[l][l]{\scriptsize DCAR}
    \psfrag{NCDR---}[l][l]{\scriptsize NCDR}
    \psfrag{sample}[t][c]{\footnotesize $N$}
    \psfrag{yaxis}[b][c]{\footnotesize $\log_{10} \delta$}
    \includegraphics[width=\textwidth]{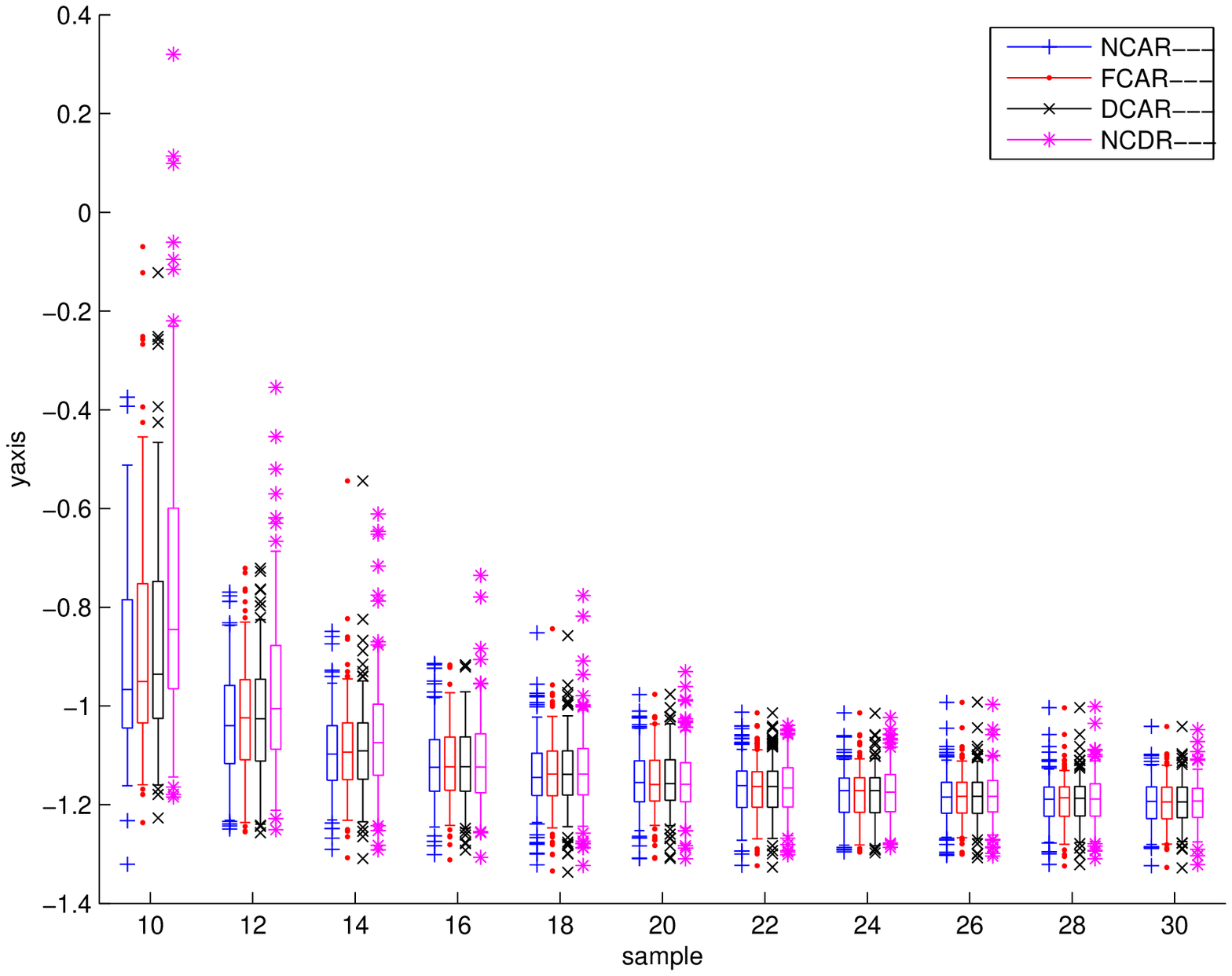}
    \caption{Magnetometer}
  \end{subfigure}
  \caption{Reconstruction error for test data for different numbers of sample
  sets. Best viewed in color.}
  \label{fig:results:sample_test}
\end{figure*}

The field components $g_z$, $h_x$ and $h_z$ are sampled uniformly over
$[-1.5,-0.5]$, $[0.5,1.5]$ and $[-1.5,1.5]$, respectively. Parameters for the
accelerometer and magnetometer are created in a similar fashion. Each component
$b_{s,i}$ of the biases has a random value sampled uniformly in $[-1,1]$. The
gain matrices are given by
\begin{equation}
\label{eq:experiment_gain}
  K_s = I_3 + \delta K_s, ~\delta K_s \sim \MCU([-0.1,0.1]),
\end{equation}
and the covariances are given by
\begin{equation}
\label{eq:experiment_covariance}
  \Sigma_s = \alpha_s
  \begin{dcases}
    \Sigma_{s,ij} \sim \MCU([0.5,2]), & \text{if } i=j \\
    \Sigma_{s,ij} \sim \MCU([-0.2,0.2]), & \text{otherwise},
  \end{dcases}
\end{equation}
where $\alpha_s = 10^\epsilon, ~\epsilon \sim \MCU([-2,-4])$, which guarantees
that the covariance matrices are positive definite, but allows high condition
numbers. Instead of using the nominal value $K_m$, the magnetometer is rotated
by $R'$ and mirrored by $M$, so that the resulting gain matrix is given by $K_m'
= M R' K_m$, showing that the algorithm described is able to handle any
configuration between the accelerometer and magnetometer, unlike \cite{Kok2012}.
The rotation $R'$ is chosen randomly \cite{Kuffner2004}, and $M =
\text{diag}(\delta_1, \delta_2, \delta_3), ~\delta_i \sim \MCU(\{-1,1\})$, is a
random diagonal matrix that indicates the mirroring of each component.

\begin{figure*}[!t]
  \centering
  \begin{subfigure}{\scale\textwidth}
    \centering
    \psfrag{NCAR---}[l][l]{\scriptsize NCAR}
    \psfrag{FCAR---}[l][l]{\scriptsize FCAR}
    \psfrag{DCAR---}[l][l]{\scriptsize DCAR}
    \psfrag{NCDR---}[l][l]{\scriptsize NCDR}
    \psfrag{stop}[t][c]{\footnotesize $\log_{10} \gamma$}
    \psfrag{yaxis}[b][c]{\footnotesize Time ($\log_{10}$ [s])}
    \includegraphics[width=\textwidth]{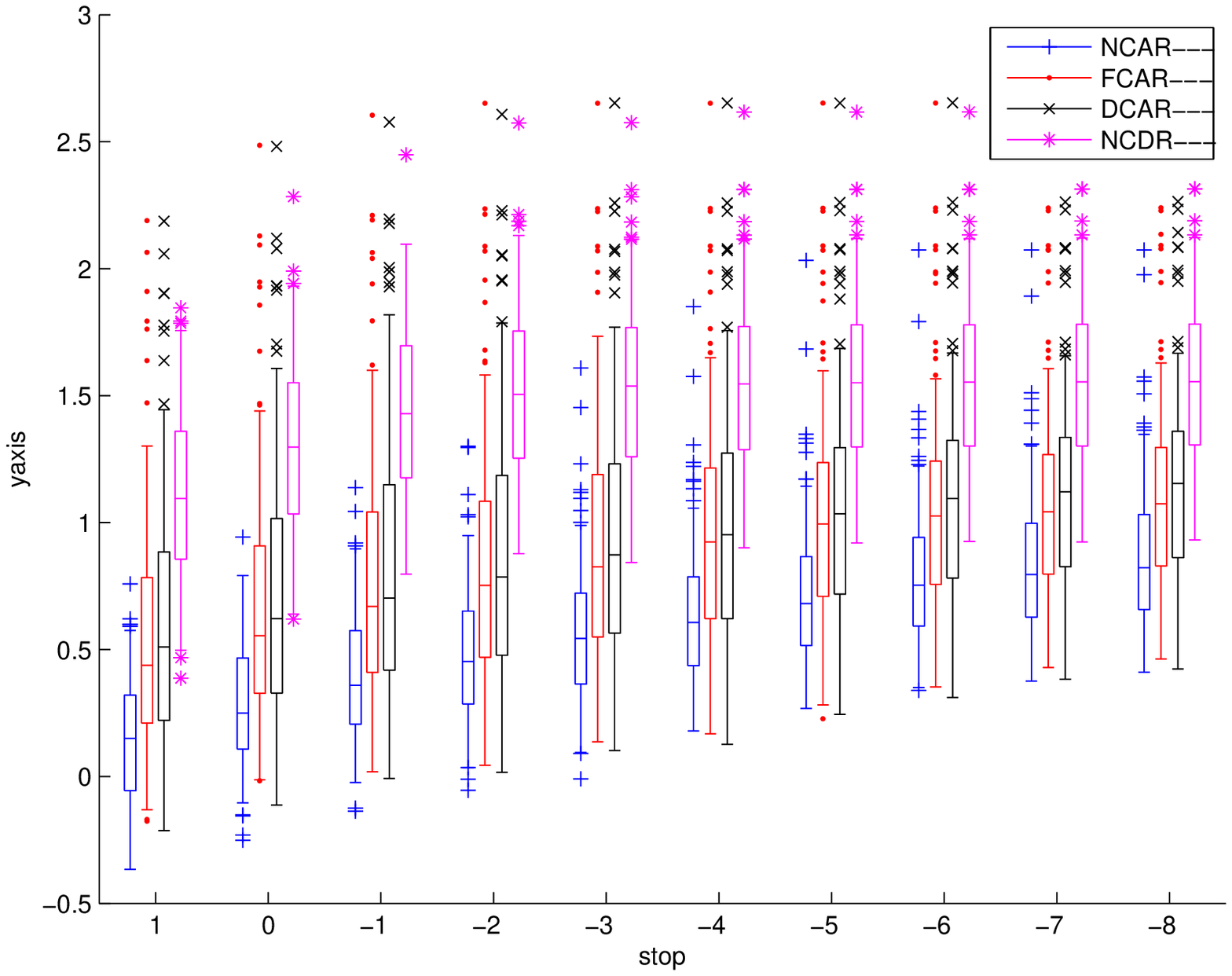}
    \caption{Varying stop condition}
    \label{fig:results:time:stop}
  \end{subfigure}
  ~
  \begin{subfigure}{\scale\textwidth}
    \centering
    \psfrag{NCAR---}[l][l]{\scriptsize NCAR}
    \psfrag{FCAR---}[l][l]{\scriptsize FCAR}
    \psfrag{DCAR---}[l][l]{\scriptsize DCAR}
    \psfrag{NCDR---}[l][l]{\scriptsize NCDR}
    \psfrag{sample}[t][c]{\footnotesize $N$}
    \psfrag{yaxis}[b][c]{\footnotesize Time ($\log_{10}$ [s])}
    \includegraphics[width=\textwidth]{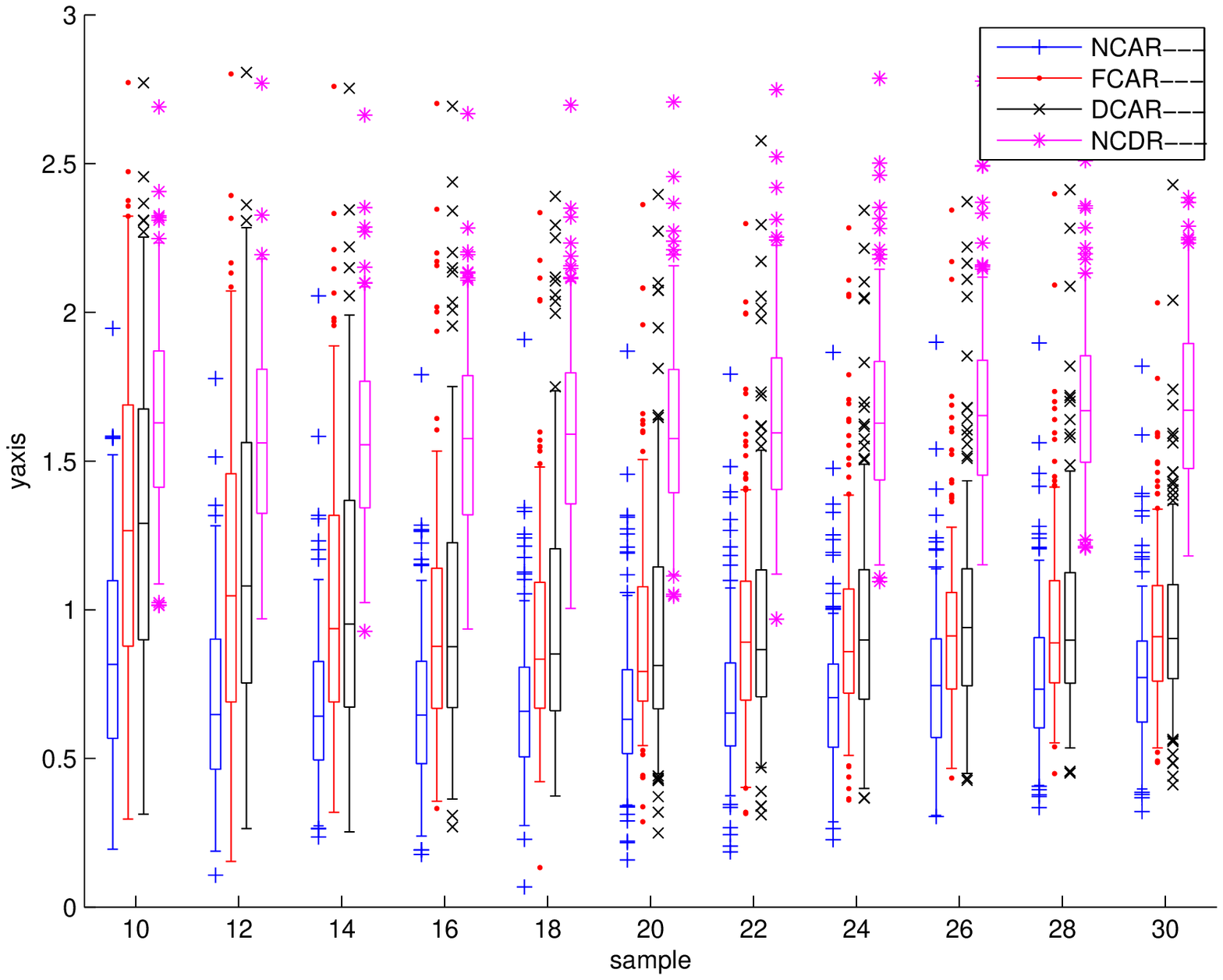}
    \caption{Varying number of samples sets}
    \label{fig:results:time:sample}
  \end{subfigure}
  \caption{Calibration time. Best viewed in color.}
  \label{fig:results:time}
\end{figure*}

After calibration, the estimated value of $\Theta$ was used to rebuild estimates
$\breve \mu_s[i]$ for $\mu_s[i]$, as in Eq.~\eqref{eq:sensor_reading_model}. To
compare the quality of the reconstructions, the following measurement was used:
\begin{equation}
\label{eq:distance}
  \delta_s = \sqrt{\frac{\sum_{i=1}^{N} \Delta[i]
  \left(\mu_s[i] - \breve \mu_s[i]\right)^T
  \Sigma_s^{-1}
  \left(\mu_s[i] - \breve \mu_s[i]\right)}{\sum_{i=1}^{N} \Delta[i]}},
\end{equation}
where $\Sigma_s$ is the real sensor covariance. This measurement is derived from
the Mahalanobis distance \cite{mahalanobis1936generalized} and represents the
average error in standard deviations.

The rotation optimization described in Sec.~\ref{sec:rotation_estimation} can be
performed in one of two ways: always using gradient descent in
Eq.~\eqref{eq:rotation_estimation} or using the approximation in
Eq.~\eqref{eq:rotation_approximation} until $J_k > J_{k-1}$ and then changing to
gradient descent in subsequent iterations. As the approximation does not
minimize the cost directly, eventually the change to gradient descent will be
made. Another variation involves re-estimating the covariance using
Eq.~\eqref{eq:covariance_estimation} or only using the initial estimate given by
Eq.~\eqref{eq:covariance_initial_estimation}. As some approaches consider the
covariance matrix to be diagonal \cite{vasconcelos2011geometric}, this
restriction is also evaluated.

Hence there are four algorithms to compare: not refitting covariance and
optimizing the rotation directly using gradient descent (NCDR); not refitting
covariance and approximating the rotation (NCAR); and fitting full covariance
and diagonal covariance and approximating the rotation (FCAR and DCAR,
respectively). The results obtained while readjusting the covariance with direct
rotation optimization are not reported as these were considerably worse.

For each sample set, another set with different samples and rotations for the
same sensor parameters was created to evaluate the generalization of the fitted
parameters, that is, if the parameters are able to fit the data not present
during calibration well. These new samples follow the description in
Sec.~\ref{sec:data_preprocessing} and provide new values $\hat \mu'[i]$, which
are used to compute the new rotations. The rotations were first approximated and
then optimized directly to reduce the cost function, as described in
Sec.~\ref{sec:rotation_estimation}.

Figure~\ref{fig:results:stop_train} shows the reconstruction error for different
stop conditions. For $\gamma < 10^{-4}$, the error does not change
significantly, justifying the choice of this value as a reference when
evaluating different number of samples. Except for NCDR, all the algorithms have
similar behavior for lower values of $\gamma$, while adjusting the covariance
was beneficial for larger stop conditions. However, as larger stop conditions
can only be justified if the computational cost is very high, it will be shown
that this difference is not relevant. NCDR has slightly worse magnetometer
calibration because of the failure of the gradient descent method to escape
local minima.

It is interesting to note the apparent magnetometer overfitting for $\gamma <
10^{-2}$, as the reconstruction error for the training set increases even though
the cost function decreases. Nonetheless, Fig.~\ref{fig:results:stop_test} shows
that the reconstruction error for the new data not used during training
generally decreases. This suggests that, although the reconstruction worsens,
the parameters improve the generalization. The significantly higher error for
one case with NCDR indicates that it is not as robust as the other methods using
approximations. For all the algorithms, about 75\% of all cases have errors of
less than 0.1 standard deviations from the correct value for both sensors,
showing that the proposed method is able to fit the parameters well.

Figure~\ref{fig:results:sample_train} shows that all the algorithms eventually
converge to solutions that are close to each other as the number of samples
increases, with errors on every run lower than 0.1 standard deviations for $N >
25$ for all methods, except for one run with NCDR. However, a diminishing return
on the number of samples is observed for $N > 20$ as the error decreases very
little if more data is collected. The NCAR algorithm appears to be the most
robust, with a lower error spread for a smaller number of samples, so that less
user data is needed to achieve acceptable performance for most situations. The
NCDR algorithm usually has the largest errors, which, as in the previous case,
are the result of poor local minima. Both algorithms with readjusted covariance
have similar performance, indicating that the diagonal simplification is a good
one. Figure~\ref{fig:results:sample_test} shows slightly better performance
overall, indicating that the estimates obtained are able to generalize and fit
new unseen samples well.

While the NCAR, FCAR and DCAR algorithms had similar performances for reasonable
stop conditions and numbers of samples, the NCDR algorithm not only had the
worst performance for all the conditions simulated but also required more
computational time, as shown in Fig.~\ref{fig:results:time:stop}. It should be
noted that, although DCAR has a simpler covariance to fit than FCAR, it has a
higher computational cost because of implementation details. Hence, the
calibration time should be about the same for these algorithms in a final
implementation. As stated earlier, using higher values of $\gamma$ and
covariance re-estimation would be justified if the computational cost was high.
However, the cost for NCAR is significantly lower, with $\gamma = 10^{-2}$
having the same cost as $\gamma=10$ for FCAR. Therefore re-estimating the
covariance using the rotation approximation does not improve performance and
instead slows the processing considerably.

\begin{figure*}[!t]
  \centering
  \begin{subfigure}{\scale\textwidth}
    \centering
    \psfrag{Samples}[t][c]{\footnotesize Samples}
    \psfrag{Values}[b][c]{\footnotesize Values}
    \includegraphics[width=\textwidth]{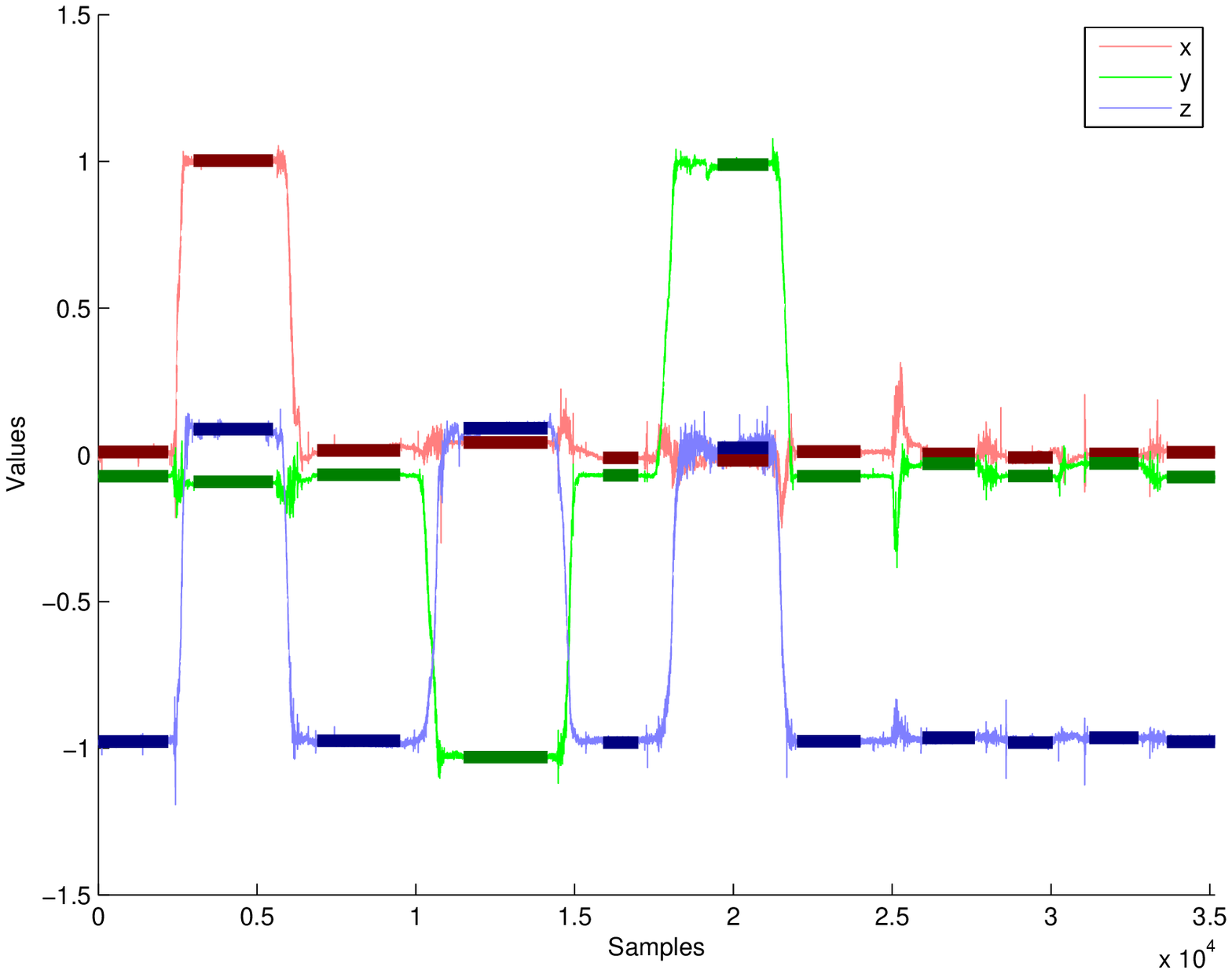}
    \caption{Accelerometer}
    \label{fig:results:real:accelerometer}
  \end{subfigure}
  ~
  \begin{subfigure}{\scale\textwidth}
    \centering
    \psfrag{Samples}[t][c]{\footnotesize Samples}
    \psfrag{Values}[b][c]{\footnotesize Values}
    \includegraphics[width=\textwidth]{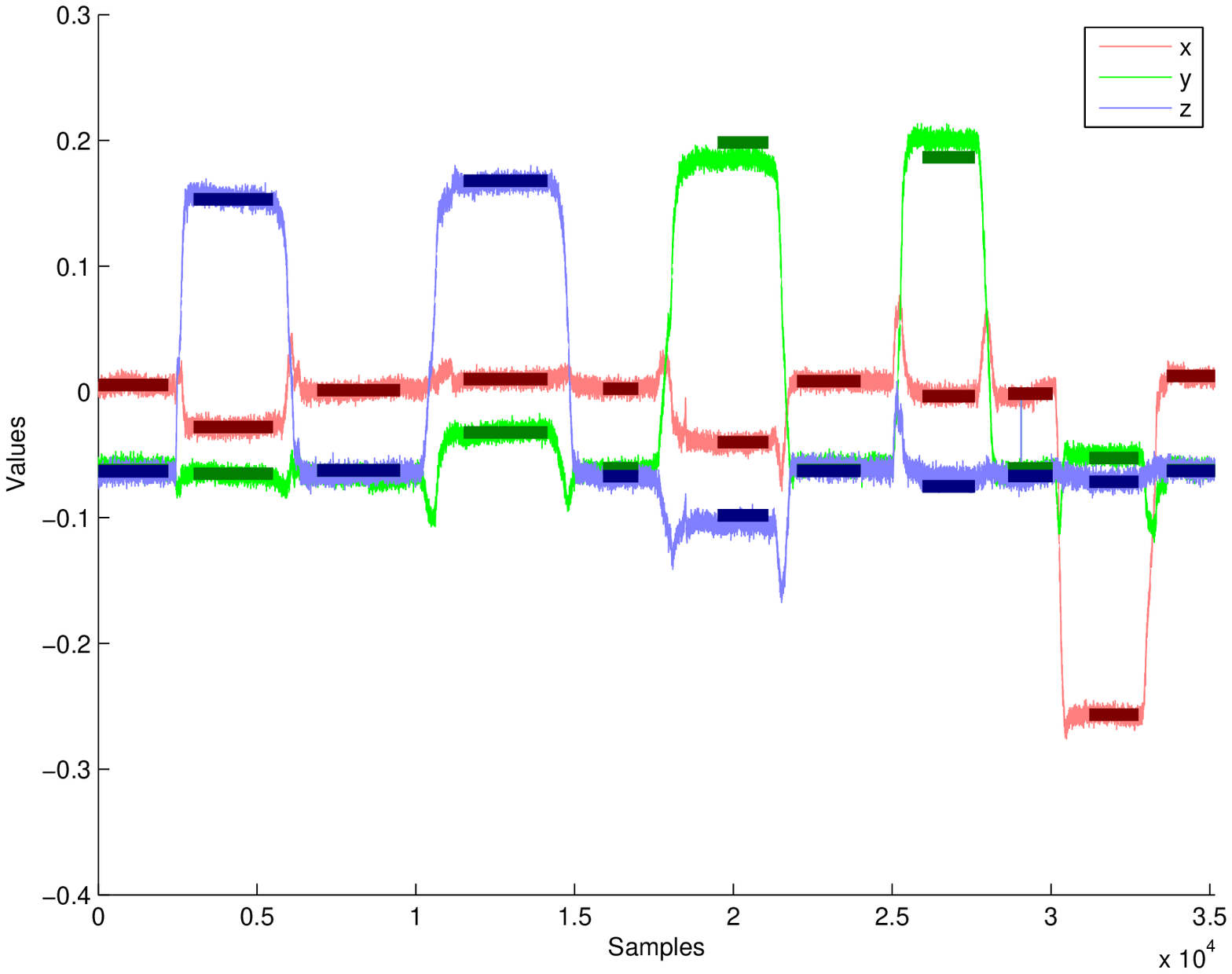}
    \caption{Magnetometer}
    \label{fig:results:real:magnetometer}
  \end{subfigure}
  \caption{Calibrated real sensors. The light lines show the values read and the
  thick, darker straight lines show the estimated reading for each set. Best
viewed in color.}
  \label{fig:results:real}
\end{figure*}

As shown in Figs.~\ref{fig:results:stop_train} and \ref{fig:results:stop_test},
$\gamma$ values of less than $10^{-4}$ do not improve the performance
significantly but cause the computing time for the approximation methods to
increase. This occurs because rotation optimization is still taking place,
which, although having little effect on the calibration, is time consuming. This
also explains why the time taken by the NCDR algorithm stabilizes.

The time for NCAR is almost constant for all the numbers of samples tested, as
shown in Fig.~\ref{fig:results:time:sample}. This occurs despite the fact that
the cost function in Eq.~\eqref{eq:partial_cost} is a sum that increases as the
number of sample sets increases and the stop condition compares two iterations
of the algorithm. As more terms are used in the sum, the difference $J_{k-1} -
J_k$ would also be expected to increase, so that it takes longer to satisfy the
stop condition. Moreover, the optimization steps also increase in complexity,
requiring more time for each step in Sec.~\ref{sec:optimization_steps}.
Nonetheless, the presence of more samples increases the speed at which the
estimates converge, requiring less iterations to achieve smaller values of $J$.
\subsection{Real calibration}
Once it has been shown that the proposed algorithm works correctly with
different simulated sensors, calibration of real sensors must be tested. The
main difference between the real and simulation experiments is that, although
the simulation allowed a wide variety of sensor parameters to be tested, it also
simulated a perfectly still system during the capture of each dataset, which was
an assumption made in Sec.~\ref{sec:data_preprocessing}. Since the NCAR
algorithm showed the best performance across all simulations, it was chosen to
calibrate the real system.

In real-world applications the system may not be perfectly still, and the
robustness of the algorithm in applications where this is the case must be
tested. Another assumption was that the magnetic field is constant throughout
the measurements, since it is assumed that it can be written as in
Eq.~\eqref{eq:fields}. This also may not be valid in the real world because of
the influence of power cables and any metal structures nearby.

To test the algorithm's robustness under these conditions, a RoboVero board was
held by hand to introduce vibration and measurements were taken inside a large
room so that magnetic disturbances were present but not overwhelming. Both the
accelerometer and the magnetometer were configured to have the same sample rate
so that both sensors took the same number of measurements for each set. The set
length was determined by using time slices in which the norms of the values read
were kept approximately constant. This also introduces errors as the underlying
values can vary without changing the norm considerably.

The results are shown in Fig.~\ref{fig:results:real}, where the darker lines are
the estimated mean sensor readings for each measurement set. These are computed
from the estimated parameters using Eq.~\eqref{eq:sensor_reading_model}, with
length equal to the number of measurements $\Delta_s[i]$ in each set. The light
lines are the real sampled values.

Fig.~\ref{fig:results:real:accelerometer} shows that the estimated values for
the accelerometer were very close to the real ones despite the sensors not being
held still during each set. This confirms the robustness of the algorithm when
the first assumption is violated and supports its use when the sensors cannot
be kept stationary.

Fig.~\ref{fig:results:real:magnetometer} shows that the magnetometer is also
calibrated correctly but that its measurements vary a lot more than the
accelerometer measurements. As a result, some measurement sets have slight
offsets, as is the case of the sets around $20000$ and $26000$ samples. These
errors are associated with distortion of the magnetic field by the building's
metallic structure. However, the small size of these offsets show that this does
not affect significantly the calibration, confirming the algorithm's robustness
in the presence of this kind of distortion.

It will be recalled that the calibration was performed using only the measured
values shown in Fig.~\ref{fig:results:real} without any external references or
sensors. Therefore, the proposed algorithm's robustness to violation of the
assumptions about the only measurements available, together with the
high-quality calibration that it provides, should make it a good candidate for
any calibration of these kinds of sensors.

\section{Conclusion}
\label{sec:conclusion}
This paper has described an algorithm to calibrate the parameters of an
accelerometer and a magnetometer using only sensor measurements without any
external information. The algorithm is designed to estimate the gain, bias, and
covariance of each sensor, as well as the orientation of each measurement and
the direction of the Earth's magnetic field. Hence, the calibration can be
performed almost anywhere by anyone as only the sensor readings are required.
This is the most generic setting for this kind of problem, since all the
parameters are calibrated and any external information can be considered a
constraint on the parameters and does not make the calibration harder.

A comparison was made of the base method and its variants, which use an
approximation of the cost function being minimized. Simulated results show that
the simplest variant is also the fastest and most robust, with the smallest
worst-case errors. All the algorithms are able to achieve an error of less than
0.1 standard deviations for the data used to calibrate the devices, indicating
that the algorithms were correctly trained, and also for new data not used in
the calibration, indicating that the parameters obtained are suitable for use
when taking sensor readings and that the learning algorithm was able to
generalize. Test with real sensors and adverse conditions that violate the
assumptions underlying the derivation of the algorithms showed that high-quality
calibration can be achieved in non-controlled settings, confirming the
robustness of the algorithm.

Future studies should investigate use of the estimated sensor rotations to
compute the rotation between body and sensor frame, as described
in~\cite{Miranda2014}, since the sensor frame used was specific to this study
and may be different from the desired frame. Another possibility for further
study is to integrate external references, such as vision systems \cite{le2009},
so that other sources of information can be used during calibration to reduce
the errors even further.

\section*{Acknowledgment}
The authors would like to thank CNPq and FAPESP for the financial support.

\bibliographystyle{templates/IEEEtran/bib/IEEEtran}
\bibliography{paper}

\begin{IEEEbiographynophoto}{Conrado S. Miranda}
  received his M.S. degree on Mechanical Engineering and his B.S. in Control
  and Automation Engineering from the University of Campinas (Unicamp), Brazil,
  in 2014 and 2011, respectively. He is currently a Ph.D. student at the School
  of Electrical and Computer Engineering, Unicamp. His main research interests
  are machine learning, multi-objective optimization, neural networks, and
  statistical models.
\end{IEEEbiographynophoto}

\begin{IEEEbiography}[{\includegraphics[width=1in,height=1.25in,clip,keepaspectratio]{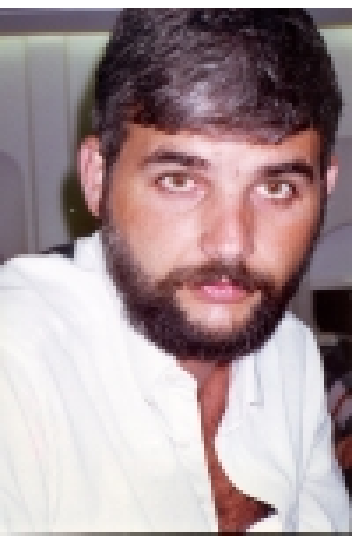}}]{Janito V. Ferreira}
  received his Ph.D. in Mechanical Engineering from the Imperial College,
  London, in 1999, and his M.S. and B.S. in Mechanical Engineering from the
  University of Campinas (Unicamp), Campinas, SP, Brazil, in 1989 and 1983,
  respectively. Since 1984 he is with the Department of Computational Mechanics
  at Unicamp, where he is presently the head of the Autonomous Mobility
  Laboratory. His main interests are in the areas of autonomous mobility,
  analysis of dynamic systems, applied and experimental mechanics, project of
  integrated system, and signal and image processing.
\end{IEEEbiography}

\end{document}